\documentclass{ifacconf}

\usepackage{graphicx}      % include this line if your document contains figures
\usepackage{natbib}        % required for bibliography
\usepackage{amsmath,amssymb,array,arydshln,color}
\usepackage{multirow}
\usepackage{diagbox}
\usepackage{subcaption}
\usepackage{caption}
\usepackage[colorinlistoftodos, textwidth=2cm, shadow]{todonotes}

\def\Ac{{\mathcal A}}

\def\Cbb{{\mathbb C}}

\def\Ec{{\mathcal E}}

\def\Gc{{\mathcal G}}

\def\Hc{{\mathcal H}}

\def\Lc{{\mathcal L}}

\def\Mc{{\mathcal M}}

\def\Nc{{\mathcal N}}

\def\Rbb{{\mathbb R}}

\def\Vc{{\mathcal V}}

\def\0{{\bf 0}}

\newcommand{\bitem}{\begin{itemize}}
\newcommand{\eitem}{\end{itemize}}
\newcommand{\btabular}{\begin{tabular}}
\newcommand{\etabular}{\end{tabular}}
\newcommand{\bcenter}{\begin{center}}
\newcommand{\ecenter}{\end{center}}
\newcommand{\bea}{\begin{eqnarray}}
\newcommand{\eea}{\end{eqnarray}}
\newcommand{\bean}{\begin{eqnarray*}}
\newcommand{\eean}{\end{eqnarray*}}
\newcommand{\ba}{\left[ \begin{array}}
\newcommand{\ea}{\\ \end{array} \right]}
\newcommand{\bear}{\begin{array}}
\newcommand{\eear}{\\ \end{array}}

\newcommand{\non}{\nonumber}

\newcommand{\ra}{\rightarrow}

\newcommand*{\QEDB}{\hfill\ensuremath{\blacksquare}}%
\newcommand{\norm}[1]{\left\lVert#1\right\rVert}

\newcounter{subequation}
\def\beasub{\addtocounter{equation}{+1}
\setcounter{subequation}{\value{equation}}
\setcounter{equation}{0}
\renewcommand{\theequation}{\arabic{subequation}\alph{equation}}
\begin{eqnarray}}
\def\eeasub{\end{eqnarray}
\setcounter{equation}{\value{subequation}}
\renewcommand{\theequation}{\arabic{equation}}}

\newtheorem{problem}{Problem}

\begin{document}
\begin{frontmatter}

\title{A Single-Adversary-Single-Detector Zero-Sum Game in Networked Control Systems\thanksref{footnoteinfo}} 

\thanks[footnoteinfo]{This work is supported by the Swedish Research Council under
	the grants 2018-04396 and 2021-06316 and by the Swedish Foundation for Strategic
	Research.}

\author{Anh Tung Nguyen,} 
\author{Andr{\'e} M. H. Teixeira,} 
\author{Alexander Medvedev}

\address{Department of Information Technology, Uppsala University, 
\\
PO Box 337, SE-751 05, Uppsala, Sweden
\\ 
(e-mail: \{anh.tung.nguyen, andre.teixeira, alexander.medvedev\}@it.uu.se).
}

\begin{abstract}                % Abstract of not more than 250 words.
This paper proposes a game-theoretic approach to address the problem of optimal sensor placement for detecting cyber-attacks in networked control systems. 
The problem is formulated as a zero-sum game with two players, namely a malicious adversary and a detector. 
Given a protected target vertex, the detector places a sensor at a single vertex to monitor the system and detect the presence of the adversary.
On the other hand, the adversary selects a single vertex through which to conduct a cyber-attack that maximally disrupts the target vertex while remaining undetected by the detector.
As our first contribution, for a given pair of attack and monitor vertices and a known target vertex, the game payoff function is defined as the output-to-output gain of the respective system.  
Then, the paper characterizes the set of feasible actions by the detector that ensures bounded values of the game payoff.
Finally, an algebraic sufficient condition is proposed to examine whether a given vertex belongs to the set of feasible monitor vertices. 
The optimal sensor placement is then determined by computing the mixed-strategy Nash equilibrium of the zero-sum game through linear programming.
The approach is illustrated via a numerical example of a 10-vertex networked control system with a given target vertex.
\end{abstract}


\end{frontmatter}
%===============================================================================

\section{Introduction}
%%----------------
The notion of networked control systems has gained popularity in modeling and analysis of real-world large-scale interconnected systems such as power systems, transportation networks, and water distribution  networks.
Networked control systems, generally employing non-proprietary and pervasive communication and information technology, such 
as the Internet and wireless communications, may leave the systems vulnerable to  cyber-attacks \citep{teixeira2015secure} and inflict significant 
financial and societal costs. Reports on Stuxnet \citep{falliere2011w32}, for example, have shown the devastating consequences of this malicious software attack on the nuclear program of Iran.
Motivated by the above observations, cyber-physical security has become an increasingly
important aspect of control systems in recent years.

\begin{figure}[!t]
	\centering
	\includegraphics[width=0.35\textwidth]{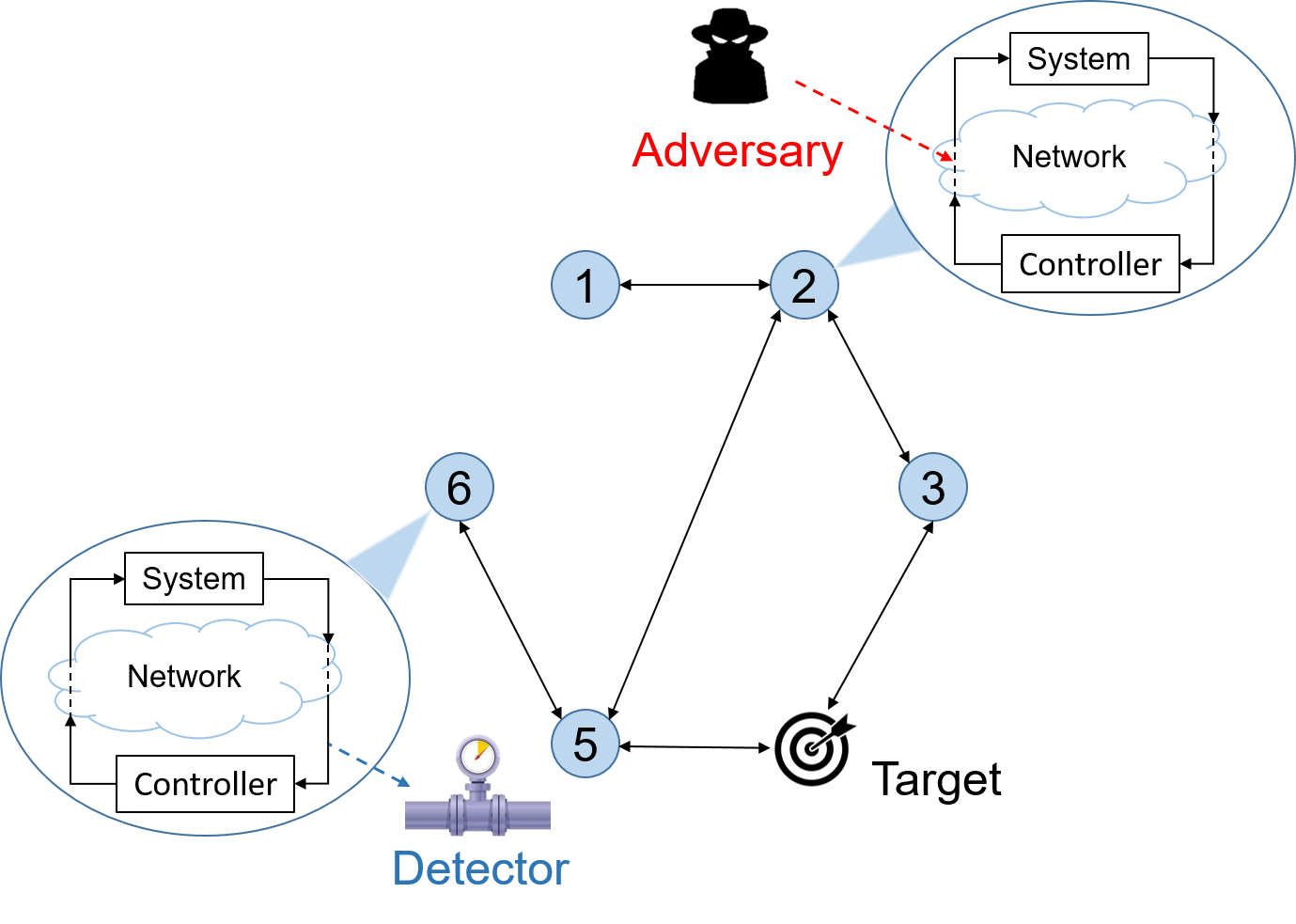}
	\caption{Visualisation of a single-adversary-single-detection zero-sum game in a networked control system.}
	\label{fig:illustration}
\end{figure}

This study considers a continuous-time networked control system under attack with two strategic agents: a malicious adversary and a detector.
The system consists of multiple one-dimensional subsystems, so-called vertices, in which there  
%are a single attack vertex, a single monitor vertex, and 
exists a single protected target vertex.
The purpose of the adversary is to affect the output of the target vertex without being detected. To this end, the adversary chooses one vertex to attack and directly injects attack signals into its input. %\textcolor{blue}{stealth?}
Meanwhile, the detector chooses one monitor vertex and measures its output, with the aim of unmasking the presence of the adversary.
Assuming both agents to be strategic, we investigate the optimal selection of the monitor vertex through a game-theoretic approach.
%We adopt a game-theoretic approach to deal with this problem.
%%
Fig.~\ref{fig:illustration} visualizes the above-defined game in a networked control system.
%%

%%----------------
The game-theoretic approach has been successfully applied to tackle the problem of robustness, security, and resilience of cyber-physical systems \citep{zhu2015game}.
%%
%\citet{zhu2015game}  
It allows us to deal with the robustness and security of cyber-physical systems within the common well-defined framework  of $\Hc_\infty$ robust control design.
Further, many other concepts of games describing networked systems subjected to cyber-attacks such as dynamic games \citep{gupta2016dynamic} and stochastic games \citep{miao2018hybrid}
%and network monitoring games \citep{milovsevic2019network} 
have been recently studied.

Although the above games were successful in studying control systems subjected to cyber-attacks such as denial-of-service attacks, changing the locations of detectors to increase the detection of such cyber-attacks was not considered.
To address this gap, \citet{pirani2021game} consider a game-theoretic formulation where the defender chooses the location of sensors in a networked system, to protect against an adversary that aims at maximally disrupting the system while remaining undetected.
The game payoff in \citet{pirani2021game} has been formulated by combining the maximum $\mathcal{L}_2$ gains of multiple outputs w.r.t. a single input representing the attack signal.
On the one hand, these multiple $\mathcal{L}_2$ gains are evaluated separately and thus may be attained for different optimal input signals, possibly resulting in pessimistic payoffs that cannot be attained by any admissible input signal.  
On the other hand, the use of a maximum gain for characterizing detectability corresponds to an optimistic perspective, where the adversary attempts to maximize the energy of the detection output, instead of the opposite.
%%
%%----------------

%\tcr{
In this paper, we consider a game-theoretic approach that is inspired and related to the one in \citet{pirani2021game}.
However, to address the above-mentioned limitations, we invoke the output-to-output gain (OOG) proposed in \citet{teixeira2015strategic,teixeira2021security} as the game payoff for the adversary and the detector.
This game payoff affords us to fully explore the cyber-attack impact on the monitor and the target outputs simultaneously with a single input signal. 
As our main contributions, we cast the optimal selection of a monitor vertex as a zero-sum game and investigate the existence of a set of feasible monitor vertices that, if selected, result in a bounded game-payoff. 
We show that the existence of such a set is related to the system-theoretic properties of the underlying dynamical system, namely its relative degrees.
Then, we propose an algebraic condition to characterize the set of feasible monitor vertices that guarantee a bounded game payoff for any attack vertex.
Finally, a numerical example is given to demonstrate the effectiveness of the proposed approach. 
Further, a mixed-strategy Nash equilibrium of the game is also investigated in a simulation example.
%}

%%--------------
We conclude this section by providing the notation to be used throughout this paper.
The problem formulation is introduced in Section \ref{sec:prob_back}.
Thereafter, Section \ref{sec:nec_suf_cond} investigates and characterizes the set of feasible monitor vertices through the system-theoretic properties of the system.
%%
%In particular, we derive an algebraic condition to seek all the feasible monitor vertices without knowing which vertex is attacked.
%%
Section~\ref{sec:num_eg} presents a numerical example of the zero-sum game between an adversary and a detector and computes the optimal monitor selection based on a mixed-strategy Nash equilibrium.
%%
%Concluding remarks are provided in 
Section \ref{sec:concl} concludes the paper.

%%--------------
{\bf Notation:} 
the set of real positive numbers is denoted as $\Rbb_+$ ; $\Rbb^n$ and $\Rbb^{n \times m}$ stand for sets of real n-dimensional vectors and n-row m-column matrices, respectively.
Let us define $e_i \in \Rbb^n$ with all zero elements except the $i$-th element is set as $1$.
%, e.g., $e_i = \big[\underbrace{0,~0,\ldots}_{i-1},1,\ldots,~0\big]^\top \in \Rbb^n$.
%%
A continuous-time system with the state-space model $\dot{x}(t) = Ax(t) + Bu(t),\; y(t) = Cx(t) + Du(t)$ is denoted as $\Sigma \triangleq (A,B,C,D)$.
Consider the norm $\norm{x}_{\Lc_2 [0,T]}^2 \triangleq \int_{0}^{T} \norm{x(t)}_2^2~dt$.
The space of square-integrable  functions is defined as $\Lc_{2} \triangleq \bigl\{f: \Rbb_{+} \rightarrow \Rbb ~|~ \norm{f}_{\Lc_2 [0,\infty]} < \infty \bigr\} $
and the extended space be defined as $\Lc_{2e} \triangleq \bigl\{ f: \Rbb_{+} \rightarrow \Rbb ~|~ \norm{f}_{\Lc_2 [0,T]} < \infty,~ \forall~ 0 < T < \infty \bigr\} $.
%%
%\\
Let $\Gc \triangleq (\Vc, \Ec, \Ac)$ be a digraph with the set of $N$ vertices $\Vc = \{v_1, v_2,...,v_N\}$,
the set of edges $\Ec \subseteq \Vc \times \Vc $, and the  adjacency matrix $A = [a_{ij}]$.
For any $(v_i,v_j) \in \Ec, ~i\neq j$, the element of the adjacency matrix $a_{ij}$ is positive, 
and with $(v_i,v_j) \notin \Ec$ or $i = j$, $a_{ij} = 0$. 
The degree of vertex $v_i$ is denoted as 
$d_i =  \sum_{j=1}^{n} a_{ij}$ and the degree matrix of graph $\Gc$ is defined as 
$D = {\bf diag}\big(d_1, d_2,\dots,d_N\big)$, where ${\bf diag}$ stands for a diagonal matrix.
The Laplacian matrix is defined as $L = [\ell_{ij}] = D - A$.
Further, $\Gc$ is called an undirected graph if  $A$ is symmetric.
An edge of an undirected graph $\Gc$ is denoted by a pair $(v_i,v_j) \in \Ec$. 
An undirected graph is connected if for any pair of vertices there exists at least one path between two vertices.
The set of all neighbours of vertex $v_i$ is denoted as $\Nc_i = \{v_j \in \Vc: (v_i,v_j) \in \Ec \}$.
%\begin{lem} \label{lemma:Lsemi}
%	\citep{oh2015survey}
%	If a graph $\Gc$ is connected, then its Laplacian matrix $L$ is positive semidefinite.
%	%%
%	Moreover, $z^T L z = 0$ if and only if $z = a {\bf 1}_n$ for some $a\in \Rbb$.
%	%% 
%	%Finally, we have $0 \leq \rho_2(\Lc) K_n \leq \Lc $, where $K_n = I_n - \frac{1}{n} {\bf 1}_n^T {\bf 1}_n$ .
%\end{lem}
%%

%------------------------------------------------
%------------------------------------------------
\section{Problem formulation} % Unnumbered section
\label{sec:prob_back}
%------------------------------------------------
%------------------------------------------------
%\AT{add one sentence explaining what this seciton is about.}
This section consists of three subsections.
Firstly,
the networked control system in the presence of a cyber-attack is defined. 
Then, we introduce the optimal stealthy data injection attack which will be studied throughout the paper.
The last subsection describes our game-theoretic approach to select feasible monitor vertices.
%------------------------------------------------
%------------------------------------------------
\subsection{Networked control system under attack}
Consider a connected undirected network $\Gc \triangleq (\Vc, \Ec, \Ac)$ with $N$ vertices, 
the state-space model of a one-dimensional vertex $v_i$ is described:
\begin{align}
	\dot x_i(t) = u_i(t), ~~ i \in \bigl\{1,~2,\ldots,~N\bigr\},
	\label{sys:xi1}
\end{align}
where $x_i(t) \in \Rbb$ is the state of vertex $v_i$. 
Due to the fact that states of all the vertices are not always available,
we employ the widely-used displacement-based control law for networked control systems: 
%in \citet{oh2015survey}:
%%
\begin{align}
	u_i(t) = \sum_{j \in \Nc_i} \big(x_j(t) - x_i(t)\big).
	\label{sys:u}
\end{align}
For convenience, let us denote $x(t)$ as the state of the networked control system, $x(t) = \big[x_1(t),~x_2(t),\ldots,~x_N(t)\big]^\top$.
In our setup, 
the adversary conducts time-dependent malicious action $a(t) \in \Rbb$ at the input of vertex $v_a$:
\begin{align}
	u_a(t) = \sum_{j \in \Nc_a} \big(x_j(t) - x_a(t)\big) + a(t).
	\label{sys:ua}
\end{align}
The purpose of the adversary is to manipulate the output of a given target vertex $v_\tau$. On the other hand, the detector places a sensor at the output of vertex $v_m$ to monitor attack signals. 
The system model \eqref{sys:xi1} under the control law \eqref{sys:u} can be rewritten in the presence of attack signals at the vertex $v_a$ \eqref{sys:ua} with two outputs observed at the two vertices $v_\tau$ and $v_m$:
\begin{align}
	\dot x(t) &= - L x(t) + e_a a(t),
	\label{sys:x}
	\\
	y_\tau (t) &= e_\tau^\top x(t),
	\label{sys:yt}
	\\
	y_m (t) &= e_m^\top x(t).
	\label{sys:ym}
\end{align}

In the scope of this study, we mainly focus on the stealthy data injection attack.
This attack will be defined as follows.
Consider the above structure of the continuous-time system \eqref{sys:x}-\eqref{sys:ym}, which we denote as $\Sigma_{\tau,m} \triangleq (-L,e_a,[e_\tau,e_m]^\top,0)$, with target output $y_\tau (t) = e_\tau^\top x (t)$ and monitor output $y_m (t) = e_m^\top x(t)$.
The input signal $a(t)$ of the system $\Sigma_{\tau,m}$ is called the stealthy data injection attack if the monitor output satisfies $\norm{y_m}_{\Lc_{2}[0,T]}^2 < \delta$, in which $\delta > 0$ is called an alarm threshold.
Further, the impact of the stealthy data injection attack is measured via the energy of the target output over the horizon $[0,T]$, i.e., $\norm{y_\tau}_{\Lc_{2}[0,T]}^2$.
Without loss of generality, let us set the alarm threshold $\delta = 1$ in the remainder of this study. 
%\textcolor{blue}{You can probably drop the squares then.}
%%

The worst-case impact of the stealthy data injection attack will be further investigated in the next subsection.

\subsection{Optimal stealthy data injection attack}
The adversary attacks vertex $v_a$ with the objective of maximizing impact on the output of the target vertex $v_\tau$ while remaining undetected at the monitor vertex $v_m$, which can be formulated as the following non-convex optimal control problem \citep{teixeira2021security}:
\begin{align}
% 	0 \leq \gamma^\star ~ \triangleq ~
	J_{\tau}(v_a,v_m) \triangleq ~
	\underset{a \in \Lc_{2e},\; x(0) = 0}{\sup} ~
	&\norm{y_\tau}_{\Lc_2}^2
	\label{opt1}
	\\
	\text{s.t.}~~
	&\norm{y_m}_{\Lc_2}^2 \leq 1.
	\non	
\end{align}
Following the details in \citet{teixeira2021security}, the above optimal control problem can be equivalently rewritten as the following optimization problem
\begin{align}
	J_{\tau}(v_a,v_m) \triangleq ~ \underset{\gamma \in \Rbb_+}{\min} &~~~~~~	\gamma
	\label{opt11}
	\\
	\text{s.t.} & ~~\norm{y_\tau}_{\Lc_2}^2 \leq \gamma \norm{y_m}_{\Lc_2}^2, ~~ \forall ~ a \in \Lc_{2e},
	\non \\ 
	& ~~~~~~~~~~~~~~~~~~~~~~~~~~~~~~
	 x(0) = 0,
	\non	
\end{align}
where the constraint may in turn be replaced with a convex Linear Matrix Inequality (\citet[Ch. 6.4]{teixeira2021security} and references therein), yielding a convex optimization problem that computes $J_{\tau}(v_a,v_m)$. 
%\textcolor{blue}{explain zero initial conditions}

\begin{rem}
	With a similar scenario,
	another objective function based on $\Lc_2$-gain for the adversary and the detector has been proposed in \citet[Sec. 3]{pirani2021game}.
	The objective function in \citet{pirani2021game} was formulated in terms of the maximal $\Lc_2$-gains from the attack to the target vertices and from the attack to monitor vertices. More specifically, the objective function in \citet{pirani2021game} is given by
	$$G_\tau(v_a, v_m) = \sup_{ \|a\|_{\Lc_2} \neq 0}
	\frac{\|y_\tau\|_{\Lc_2}^2}{\|a\|_{\Lc_2}^2}
	 - \lambda \sup_{ \|a\|_{\Lc_2} \neq 0}
	\frac{\|y_m\|_{\Lc_2}^2}{\|a\|_{\Lc_2}^2},~ (\lambda \geq 0). $$
	The above objective in \citet{pirani2021game} also considers two different outputs $y_\tau(t)$ and $y_m(t)$, but note that the output energies are maximized separately, thus leading to two different optimal input signals $a(t)$ in the general case.
	%%
	%Thus, these two worst-case outputs rarely occur with the same input signal.
	%%
	By contrast, our objective function \eqref{opt11} investigates the worst-case attack impact that is simultaneously characterized by the two outputs $y_\tau(t)$ and $y_m(t)$ w.r.t. a single input signal $a(t)$.
\end{rem}

Next, we tackle the problem of the optimal selection of a monitor vertex through a game-theoretic approach.

%%
%%
% The worst-case impact of the stealthy data injection attack will be further investigated in the next subsection.

%%---------------------------
% \subsection{Single-adversary-single-detector zero-sum game \\ description}
\subsection{Game-theoretic approach to monitor vertex selection}
To defend against adversaries, we consider that the detector tackles the following problem.

%\AT{need to redefine the Problem environment to start with Problem 1. Maybe define a new counter for the prob environment?}
%%
\begin{problem} \label{prob:opt}
	(Optimal monitor selection) Given a target vertex and an arbitrary attack vertex, select a monitor vertex that minimizes the worst-case impact of the stealthy data injection attack at the attack vertex. %\textcolor{blue}{You appear to talk about measured impact and not the actual one}
\end{problem}

As the attack vertex is arbitrary, we formulate Problem~\ref{prob:opt} as a game between the detector and adversary, where the players choose $v_a$ and $v_m$ to respectively maximize and minimize the function $J_\tau (v_a,v_m)$ described in \eqref{opt11}.
Hence, Problem~\ref{prob:opt} is formalized as a zero-sum game with $J_\tau (v_a,v_m)$ as the game-payoff, namely
%We consider the following zero-sum game and 
%%
\begin{align}
	\underset{v_m \neq v_\tau \in \Vc}{\min}
	~ \underset{v_a \neq v_\tau \in \Vc}{\max}
	~~
	J_\tau (v_a,v_m).
	\label{prob:opt_form}
\end{align}

% Due to an arbitrary attack vertex $v_a$, the worst scenario for the detector is that the optimisation problem \eqref{prob:opt_form} admits no solution, i.e., the adversary indirectly attacks $v_\tau$ while almost leaving no trace at $v_m$.

While Problem \ref{prob:opt} investigates an optimal selection of the monitor vertex, there is no \textit{a priori} guarantee that a suitable monitor vertex exists for which \eqref{prob:opt_form} is bounded from above.
The following Problem raises a question of finding feasible monitor vertices such that the worst-case impact of the stealthy data injection attack is bounded.
\begin{problem} \label{prob:fea}
	(Feasible monitor vertices) Given a target vertex and an arbitrary attack vertex, find a set of feasible monitor vertices such that the worst-case impact of the stealthy data injection attack is bounded.
\end{problem}

% \begin{defn} 
% 	\label{def:set_monitor}
% 	(Set of candidate monitor vertices)
% 	Given a target vertex, a vertex is called a candidate monitor if its selection as a monitor vertex guarantees the feasibility of the optimisation \eqref{opt11} for an arbitrary attack vertex.
% \end{defn}

Formally, a set of feasible monitor vertices w.r.t. the target vertex $v_\tau$ is defined as $\Mc_\tau \triangleq \{v_m \neq v_\tau \in \Vc ~|~
	\underset{v_a \neq v_\tau \in \Vc}{\max}
	~
	J_\tau (v_a,v_m) < \infty
	\}$.
%%
% Due to an arbitrary attack vertex $v_a$, the worst scenario for the detector is that the optimisation problem \eqref{prob:opt_form} admits no solution, i.e., the adversary indirectly attacks $v_\tau$ while almost leaving no trace at $v_m$.
%%

By definition, if $\Mc_\tau \neq \emptyset$, the detector may select a vertex $v_m \in \Mc_\tau$ to ensure that \eqref{opt11} is feasible for any attack vertex, which in turn guarantees that the zero-sum game \eqref{prob:opt_form} admits a bounded value. Furthermore, characterizing the set $\Mc_\tau$ allows us to restrict the possible choices of the detector to $\Mc_\tau$. Hence, by addressing Problem~\ref{prob:fea} and characterizing $\Mc_\tau$, we can tackle Problem~\ref{prob:opt} by reformulating the zero-sum game~\eqref{prob:opt_form} as
\begin{align}
	\underset{v_m \neq v_\tau \in \Vc}{\min}
	~\underset{v_a \neq v_\tau \in \Vc}{\max}
	&~~
	J_\tau (v_a,v_m) % \tcb{~ < \infty}
	\label{prob:opt_form1}
	\\
	\text{s.t.}
	&~~~
	v_m \in \Mc_\tau. %\tcb{\neq \emptyset} ~ .
	\non
\end{align}
%%
% The above optimisation problem \eqref{prob:opt_form1} admits a finite solution that is also a solution to Problem \ref{prob:opt} if and only if the set $\Mc_\tau$ is non-empty.
% How to find elements of $\Mc_\tau$ is the question raised in Problem \ref{prob:fea}.
%%
%
The next section characterizes the set of feasible monitor vertices $\Mc_\tau$ by investigating the feasibility of~\eqref{opt11} with respect to system-theoretic properties of the dynamical system \eqref{sys:x}-\eqref{sys:ym}.
%%

%------------------------------------------------
%------------------------------------------------
%\subsection{Necessary and Sufficient condition for $\Mc_\tau$}
\section{Characterization of feasible monitor vertices}
\label{sec:nec_suf_cond}
%------------------------------------------------
%------------------------------------------------
%\AT{add some text here, explaining what we do in this section.}
% \tcb{
% This section will introduce a necessary and sufficient condition for characterising possible candidate monitor vertices.
% }
%%
Let us denote the continuous-time systems 
$\Sigma_\tau \triangleq (-L,e_a,\allowbreak e^\top_\tau,0)$ and $\Sigma_m \triangleq (-L,e_a,e^\top_m,0)$. 
%\textcolor{blue}{Cannot see where this $\Sigma$ notation is defined.}
Inspired by \citet[Th. 2]{teixeira2015strategic}, the feasibility of the optimization problem \eqref{opt11} is related to the invariant zeros of $\Sigma_\tau$ and $\Sigma_m$, which are defined as follows.

\begin{defn}  \label{def:invariant_zero}
	(Invariant zeros)
	%\citep{emami1982computation} 
	Consider the strictly proper system $\Sigma \triangleq (A,B,C,0)$ with $A,B,$ and $C$ are real matrices with appropriate dimensions. A tuple $(\lambda,\bar{x},g) \in \Cbb \times \Rbb^N \times \Rbb$ is a zero dynamics of $\Sigma$ if it satisfies
	\begin{align}
		\ba{cc}
		\lambda I - A ~~~~ & -B \\
		C & 0
		\ea
		\ba{c}
		\bar{x} \\ g
		\ea
		=
		\ba{c}
		0 \\ 0
		\ea,
		~~~ \bar{x} \neq 0.
		\label{def:inv_zero}
	\end{align}
	In this case, a finite $\lambda$ is called a finite invariant zero of $\Sigma$.
	Further, the strictly proper system $\Sigma$ always has at least one invariant zero at infinity \citep[Ch. 3]{franklin2002feedback}. %\textcolor{blue}{zero at infinity}
\end{defn}
More specifically, the optimization \eqref{opt11} is feasible if and only if the unstable invariant zeros of $\Sigma_m$ are also invariant zeros of $\Sigma_\tau$ \citep[Th. 2]{teixeira2015strategic}.
To derive a necessary and sufficient condition characterizing the set of feasible monitor vertices, we will investigate both finite and infinite invariant zeros of the two systems $\Sigma_m$ and $\Sigma_\tau$.
The following Lemma considers the former.
\begin{lem}  \label{lemma:no_un_zero}
	Consider a networked control system associated with a connected undirected graph $\Gc \triangleq (\Vc,\Ec,\Ac)$, whose vertex dynamics and control law are described in \eqref{sys:xi1} and \eqref{sys:u}, respectively.
	Suppose that the networked control system is driven by the stealthy data injection attack \eqref{sys:ua} at a single attack vertex $v_a$, and observed by a single monitor vertex $v_m$, resulting in the state-space model $\Sigma_m \triangleq (-L,e_a,e^\top_m,0)$.
	Then, the networked control system $\Sigma_m$ has no 
	invariant zero on the closed positive real line.
\end{lem}
%%
% \tcr{A proof of Lemma \ref{lemma:no_un_zero} is omitted due to limited space.} 
%\begin{pf}
%	Invariant zeros of the system $\Sigma_m = (-L,e_a,e^\top_m,0)$ have two possible cases, i.e., non-eigenvalues and eigenvalues of negative Laplacian matrix $-L$.
%	%%
%	In the first case, since invariant zeros never equal to zero \citep[Theorem 3.7]{briegel2011zeros}, the other eigenvalues of negative Laplacian matrix $-L$ might be invariant zeros of the system $\Sigma_m$.
%	%%
%	According to Lemma \ref{lemma:Lsemi}, such eigenvalues are negative.
%	%%
%	In the second case, If the invariant zeros are not eigenvalues of negative Laplacian matrix $-L$, they have negative real parts \citep[Theorem 3.9]{briegel2011zeros}.
%	%%
%	Hence, the system $\Sigma_m$ has no finite invariant zero on open right half-plane $\Rbb_{+}$.
%\end{pf}
%%
\begin{pf}
The proof follows directly from the results in \citet[Th. 3.7 \& 3.9]{briegel2011zeros}.
\QEDB
\end{pf}

%\AT{The proof in the cited paper only states there are no zeros on the real line. We need to revise this for the complex plane, or to add self-loops.}
%%
\begin{rem}
\label{rem:no_unstable_zero}
     By inheriting the results in \cite{torres2015graph}, the other invariant zeros on the closed right half-plane possibly exist if the input and output vertices have short weak paths and long strong paths between them. On the other hand, the graph representing the system \eqref{sys:xi1} under control law \eqref{sys:u} is unweighted, i.e., all the edges have the same weight values, conflicting with the above sufficient condition in \cite{torres2015graph}. %%
     We leave the necessary condition under which the system $\Sigma_m$ has no invariant zeros on the closed right half-plane for the future research.
\end{rem}
Based on the above remarks, assuming that the system $\Sigma_m$ has no finite unstable zero, we then investigate the infinite zeros of the systems $\Sigma_m$ and $\Sigma_\tau$. In the investigation, we make use of known results connecting infinite invariant zeros mentioned in Definition \ref{def:invariant_zero} and the relative degree of a linear system, which is defined below.
\begin{defn} \label{def:rela_deg}
	(Relative degree)
	\citep[Ch. 13]{khalil2002nonlinear} Consider the strictly proper system $\Sigma \triangleq (A,B,C,0)$ with $A \in \Rbb^{n \times n}$, $B$, and $C$ are real matrices with appropriate dimensions.
	The system $\Sigma$ is said to have relative degree $r ~ (1 \leq r \leq n) $ if the following conditions satisfy
	\begin{align}
		&C A^{k} B = 0, ~~ 0 \leq k < r-1,
		\non \\
		&C A^{r-1} B \neq 0.
		\label{def_red}
	\end{align}
\end{defn}
%%
%\AT{add a known result relating the infinite zeros to the relative degree, or derive the result as a lemma.}
%%
\begin{rem} \label{rem:re_deg_tf}
    Let $H(s) = C(sI-A)^{-1}B$ be the transfer function of the above system $\Sigma$.
The relative degree $r$ of the system $\Sigma$ defined in Definition \ref{def:rela_deg} is also the difference between the degrees of the denominator and the numerator of $H(s)$ \citep{khalil2002nonlinear}, which in turn corresponds to the degree of the infinite zero~\cite[Ch. 3]{franklin2002feedback}.  
\end{rem}

Based on Definition~\ref{def:rela_deg}, let us denote $r_{\tau a}$ and $r_{ma}$ as the relative degrees of $\Sigma_\tau$ and $\Sigma_m$, respectively.
%from the attack vertex $v_a$ to the target vertex $v_\tau$ and the monitor vertex $v_m$, respectively.
%%
%Note that $r_{\tau a}$ and $r_{m a}$ satisfy \eqref{def_red} for $\Sigma_\tau$ and $\Sigma_m$, respectively.
%%
In the scope of this study, we have assumed that the cyber-attack \eqref{sys:ua} has no direct impact on the outputs \eqref{sys:yt} and \eqref{sys:ym}, resulting in strictly proper systems
$\Sigma_\tau$ and $\Sigma_m$.
This implies that the relative degrees $r_{\tau a}$ and $r_{ma}$ of $\Sigma_\tau$ and $\Sigma_m$ are positive, yielding their infinite zeros.
Those infinite zeros will be considered to present
the following Theorem that gives us a necessary and sufficient condition for finding feasible monitor vertices $\Mc_\tau$.
\begin{thm} \label{lemma:rela_deg}
	Consider the strictly proper systems $\Sigma_\tau \triangleq (-L,e_a,e^\top_\tau,0)$ and $\Sigma_m \triangleq (-L,e_a,e^\top_m,0)$, in which the two systems have the same stealthy data injection attack input at a single attack vertex $v_a$ but different output vertices, i.e., $v_\tau$ for $\Sigma_\tau$ and $v_m$ for $\Sigma_m$.
	Suppose the systems $\Sigma_\tau$ and $\Sigma_m$ have relative degrees $r_{\tau a}$ and $r_{ma}$, respectively.
	Then, the worst-case impact of the stealthy data injection attack \eqref{opt11} is bounded
	if and only if the following condition holds
	\begin{align}
		r_{ma} \leq r_{\tau a}. \label{cond_red}
	\end{align}
\end{thm}
%%
% \tcr{
% A proof of Lemma \ref{lemma:rela_deg} is omitted due to limited space.
% }
\begin{pf}
	Followed by \citet[Th. 2]{teixeira2015strategic},
	the optimization problem \eqref{opt11} is feasible if and only if  $\Sigma_m$ has  unstable invariant zeros that are also invariant zeros of  $\Sigma_\tau$.
	Based on Lemma \ref{lemma:no_un_zero} and Remark~\ref{rem:no_unstable_zero}, $\Sigma_m$ has no finite unstable invariant zero, which leaves us to analyze infinite zeros of those systems.
    Recall the equivalence between the relative degree of a SISO system and the degree of its infinite zero. 
 	Hence, a necessary condition to guarantee the feasibility of the optimization \eqref{opt11} is that the number of infinite invariant zeros of $\Sigma_m$ is not greater than that of $\Sigma_\tau$.
 	This implies  $r_{m a} \leq r_{\tau a}$. 
 	For sufficiency, it remains to show that 
 	%%
 	%\tcr{
 	if $r_{m a} \leq r_{\tau a}$, any infinite zeros of $\Sigma_{m}$ are also infinite zeros of $\Sigma_{\tau}$.
 	We will investigate each infinite zero of $\Sigma_{m}$
    by starting from
	their transfer functions with zero initial states
	\begin{align}
		G_{\tau a}(s) &= e_\tau^\top (sI + L)^{-1} e_a 
		= \frac{P_{\tau a}(s)}{Q(s)}
		,\non \\
		G_{ma}(s) &= e_m^\top (sI + L)^{-1} e_a 
		= \frac{P_{m a}(s)}{Q(s)}
		,
		\label{trans_fcn}
	\end{align}	
	where $s \in \Cbb$.
	%%
	%Since the state-space models $\Sigma_\tau$ and $\Sigma_m$ are minimal realisation, 
	Based on Remark \ref{rem:re_deg_tf} and the minimal realisations $\Sigma_\tau$ and $\Sigma_m$, 
	$P_{\tau a}(s)$, $P_{ma}(s)$, and $Q(s)$ are the polynomials of  degrees $N-r_{\tau a}$, $N-r_{ma}$, and $N$, respectively.
	Let us denote $z_k = \sigma_k + j \omega_k \in \Cbb,~ k \in \{1,2,\ldots,r_{ma} \}$ with infinite module as infinite zeros of $\Sigma_m$. 
	Indeed, \citep{morris2010invariant} $z_k~(1 \leq k \leq r_{ma})$ is an infinite zero of maximal degree $r_{ma}$ of $\Sigma_m$ if it satisfies
	\begin{align}
	    &\lim_{\norm{z_k}_2 \ra \infty} z_k^q G_{ma}(z_k) = 0, ~ (0 \leq q \leq r_{ma} - 1), 
	    \non \\
	    &\lim_{\norm{z_k}_2 \ra \infty} z_k^{r_{ma}} G_{ma}(z_k) \neq 0.
	\end{align}
	Further, with $0 \leq q \leq r_{m a} - 1$, we also basically have
	\begin{align}
	    &\lim_{\norm{z_k}_2 \ra \infty} z_k^q G_{\tau a}(z_k) = 
        \lim_{\norm{z_k}_2 \ra \infty}  \frac{z_k^q P_{\tau a}(z_k)}{Q(z_k)} = 0. \label{Pta_inf_zero_1}
	\end{align}
	%%
% 	\begin{align}
% 	    &\lim_{\norm{s_1}_2 \ra \infty} G_{m a}(s_1) = 0,
% 	    ~
% 	    &\lim_{\norm{s_1}_2 \ra \infty} G_{\tau a}(s_1) = 0,
% 	\end{align}
% 	which implies that $s_1$ is also an infinite zero of $\Sigma_\tau$.
% 	%%
% 	Intuitively, in order to investigate the other $i$-th infinite zero $s_i$ of $\Sigma_m,~(1 \leq i \leq r_{ma})$, one has
% 	\begin{align}
% 	    &\lim_{\norm{s_i}_2 \ra \infty} \prod_{j=1}^{i-1\leq r_{ma}-1} (s_i - s_j)~ G_{ma}(s_i) 
% 	    \non \\
% 	    &= 
% 	    \lim_{\norm{s_i}_2 \ra \infty} \frac{\prod_{j=1}^{i-1\leq r_{ma}-1} (s_i - s_j) ~ P_{m a}(s_i)}{Q(s_i)} = 0,
% 	    \\
% 	    &\lim_{\norm{s_i}_2 \ra \infty} \prod_{j=1}^{i-1\leq r_{ma}-1} (s_i - s_j)~ G_{\tau a}(s_i) 
% 	    \non \\
% 	    &= 
% 	    \lim_{\norm{s_i}_2 \ra \infty} \frac{\prod_{j=1}^{i-1\leq r_{ma}-1} (s_i - s_j) ~ P_{\tau a}(s_i)}{Q(s_i)} = 0.
% 	    \label{Pta_inf_zero_1}
% 	\end{align}
	%%
	The above limit \eqref{Pta_inf_zero_1} holds because the denominator $z^{q}P_{\tau a}(s)$ is the polynomial of degree $N-r_{\tau a}+q \leq N-1 < N$, the degree of the polynomial $Q(z_k)$.
	This implies that any infinite zeros $z_k$ of maximal degree $r_{ma}$ of $\Sigma_m$ are also infinite zeros of degree $r_{ma}$ of $\Sigma_\tau$.
	\QEDB
%}
\end{pf}

The following Lemma introduces a sufficient condition under which feasible monitor vertices guarantee the feasibility of the optimization \eqref{opt11} for an arbitrary attack vertex.
\begin{lem} \label{theo:monitor}
	Consider a networked system associated with a connected undirected graph $\Gc \triangleq (\Vc,\Ec,\Ac)$, whose vertex dynamics and control law are described in \eqref{sys:xi1} and \eqref{sys:u}, respectively.
	The networked system is driven by the stealthy data injection attack at an arbitrary attack vertex $v_a$ with its impact measured at a given target vertex $v_\tau$.
	Suppose that there exists a feasible monitor vertex $v_m$ that directly connects to all the neighbors of the target vertex $v_\tau$. Then, the optimisation \eqref{opt11} is feasible. Furthermore, defining $A_m \triangleq I + A$, this vertex $v_m$ satisfies
	\begin{align}
	    e_\tau^\top A A_m e_m = e_\tau^\top A^2 e_\tau.
 		\label{opt2:cond}
	\end{align}
\end{lem}
\begin{pf}
Recalling that the relative degrees of $\Sigma_\tau$ and $\Sigma_m$ are related to the length of the shortest paths (\textit{i.e.}, distance) from $v_a$ to $v_\tau$ and $v_m$, respectively \cite[Th. 3.2]{briegel2011zeros}, the first part of the proof follows directly from the fact that a vertex $v_m$ is connected to all the neighbors of $v_\tau$. This implies
that the distance from any arbitrary attack vertex $v_a \neq v_\tau$ to $v_m$ will be less or equal to the distance from $v_a$ to $v_\tau$. Thus, the vertex $v_m$ satisfies~\eqref{cond_red}. The remainder of the proof expresses the relation between $v_m$ and $v_\tau$ in terms of the adjacency matrix and is omitted due to space limitations.
\QEDB
\end{pf}

\begin{rem}
     To seek a set of feasible monitor vertices, 
     % that guarantee the feasibility of the optimization \eqref{opt11}, 
     the algebraic condition \eqref{opt2:cond} is simply tested with all the vertices $v_m \neq v_\tau \in \Vc$.
\end{rem}

%%
% \begin{rem} \label{remark:re_dis}
% 		Based on the definition of a walk between two vertices in a graph \citep[Chapter 1.6]{bondy1976graph}, $r_{\tau a}(r_{tma})-1$  is the minimal number of edges connecting $v_a$ to $v_\tau (v_m)$.
% 		%%
% 		%Surprisingly, Eq. \eqref{re_deg} also represents Markov parameters of the systems \eqref{sys:x} \eqref{sys:yt} and \eqref{sys:x} \eqref{sys:ym} \citep[Definition 9.2.2]{datta2004numerical}.
% 		%%
% 		Alternatively, $r_{\tau a}(r_{tma})-1$ also represents the length of the shortest path between those vertices $v_a$ and $v_\tau (v_m)$ \citep[Theorem 3.2]{briegel2011zeros}.
% \end{rem}

%------------------------------------------------
%------------------------------------------------
\section{Numerical examples}
\label{sec:num_eg}
\begin{figure} [!t]
	\centering
	\includegraphics[width=0.3\textwidth]{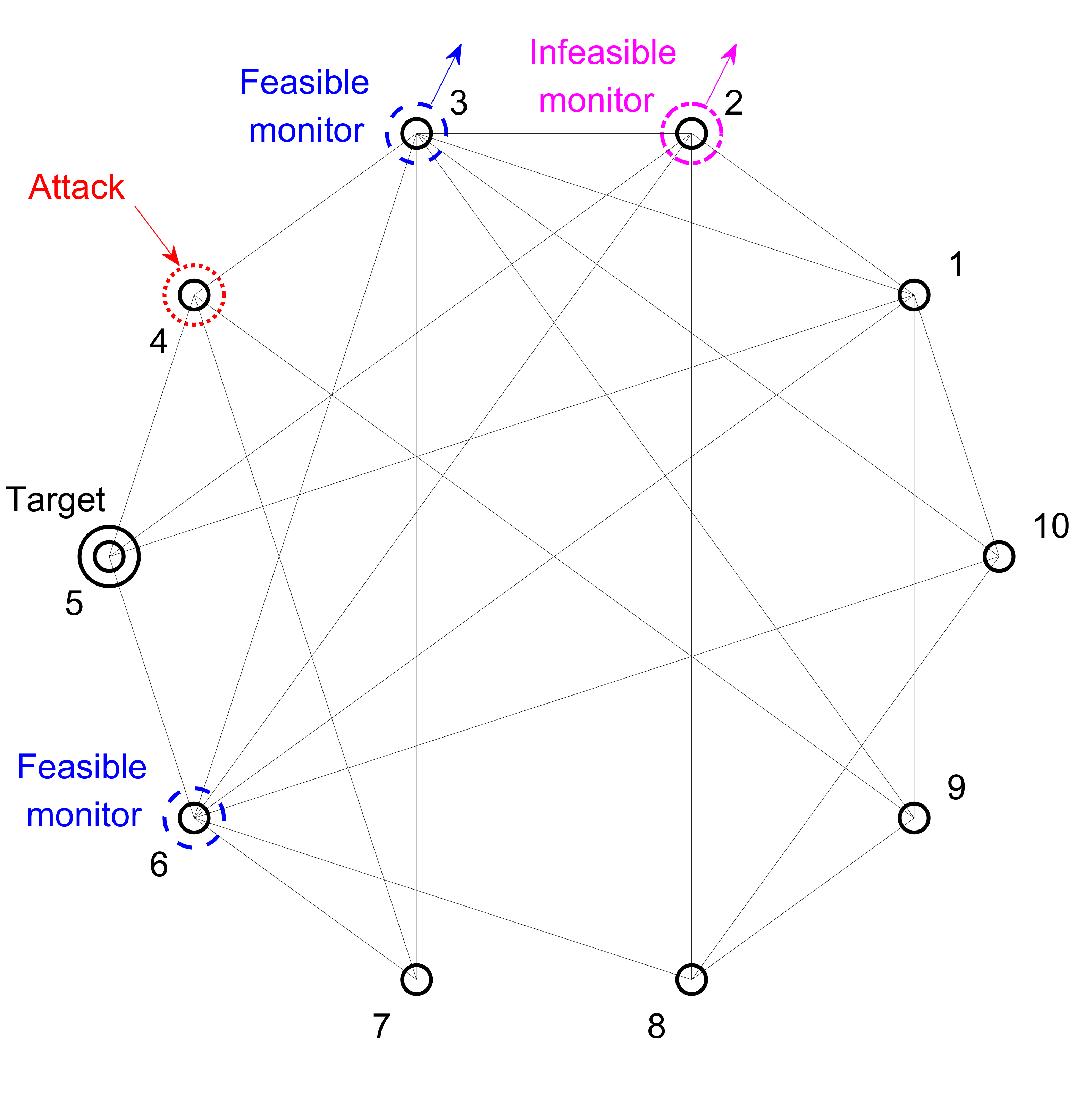}
	\caption{10-vertex networked control system with target vertex $v_5$, attack vertex $v_4$, and feasible monitor vertices $v_3$ and $v_6$.}
	\label{fig:graph_10vertex}
\end{figure}

\begin{figure} [!t]
	\begin{center}
		\begin{subfigure}{.2\textwidth}
			\centering
			\includegraphics[width=\textwidth]{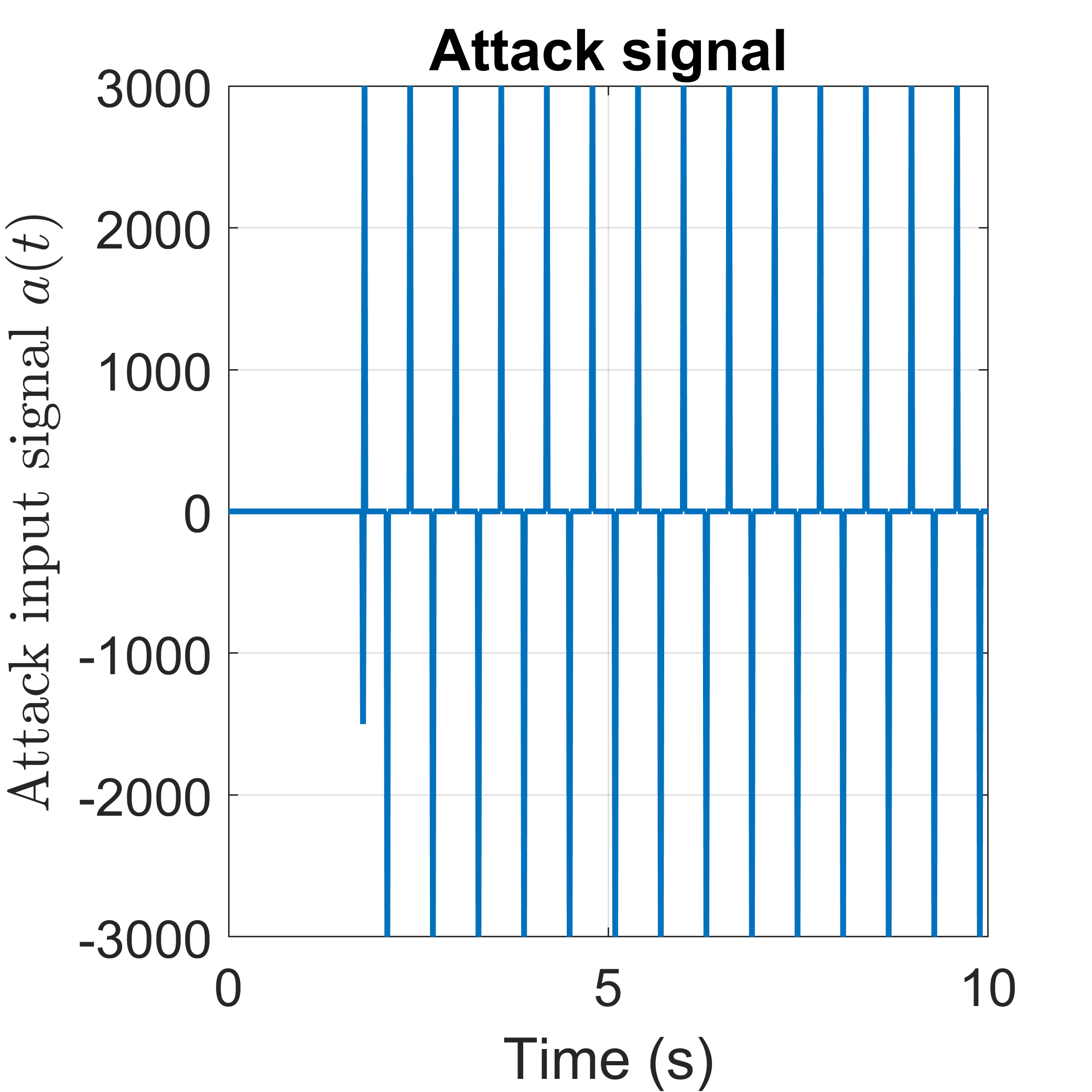}
			\caption{}
			\label{fig:attack}
		\end{subfigure}
		\begin{subfigure}{.2\textwidth}
			\centering
			\includegraphics[width=\textwidth]{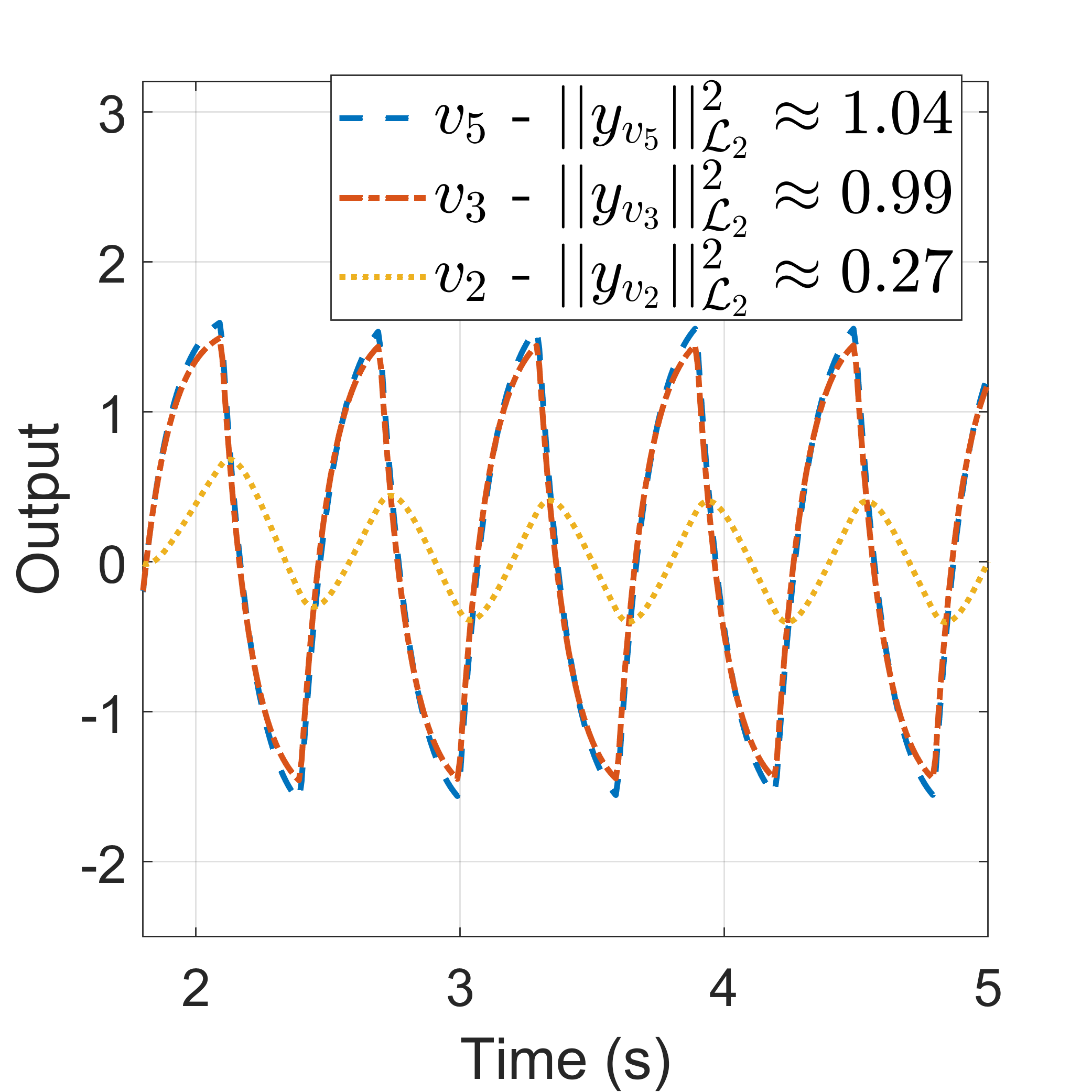}
			\caption{}
			\label{fig:output}
		\end{subfigure}
	\end{center}	
	\caption{(a) Attack signals directly injected into the input of  $v_4$ at frequency 1.67 Hz; (b) Outputs of target $v_5$, feasible $v_3$, infeasible vertices $v_2$.}
\end{figure}

\begin{figure} [!t]
    \centering
    \includegraphics[width=0.45\textwidth]{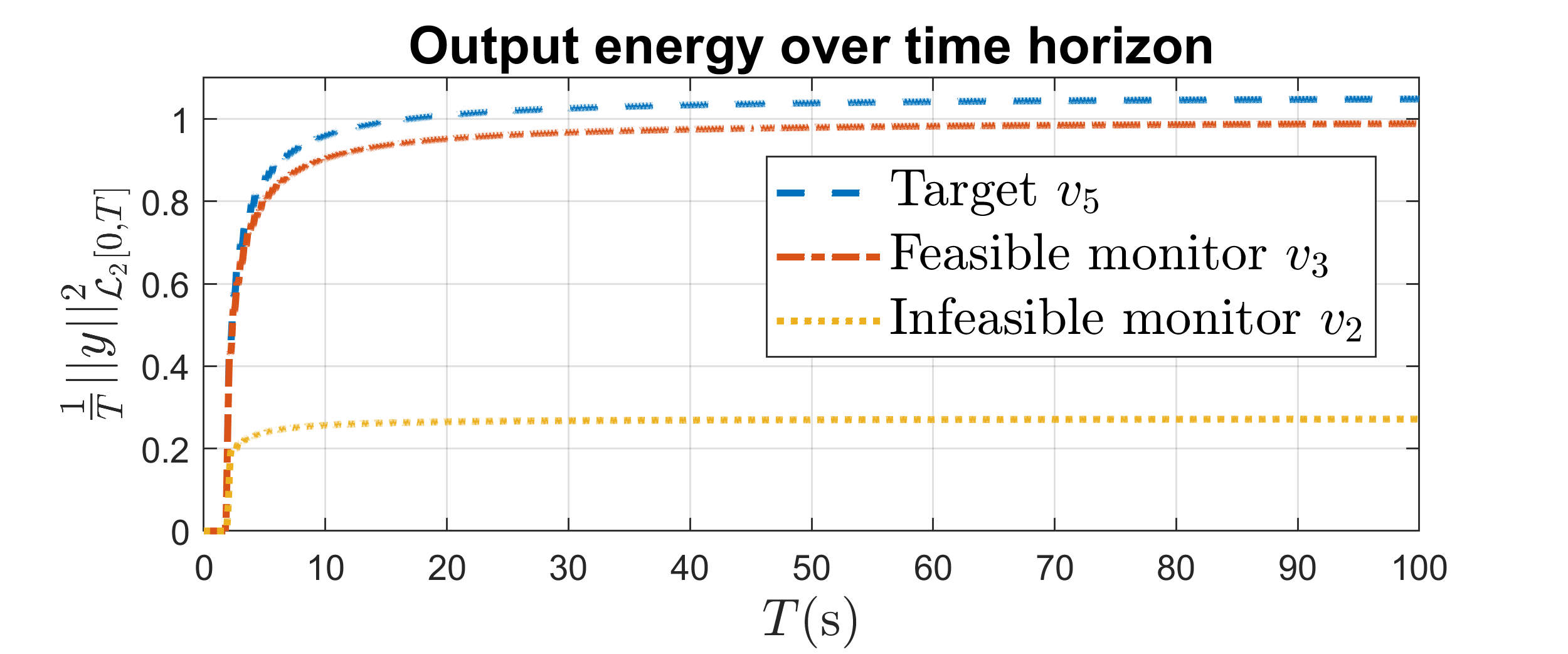}
    \caption{The output energies of vertices $v_2,~v_3,$ and $v_5$ over time horizon.}
    \label{fig:output_energy}
\end{figure}

\begin{figure}[!t]
	\centering
	\includegraphics[width=0.45\textwidth]{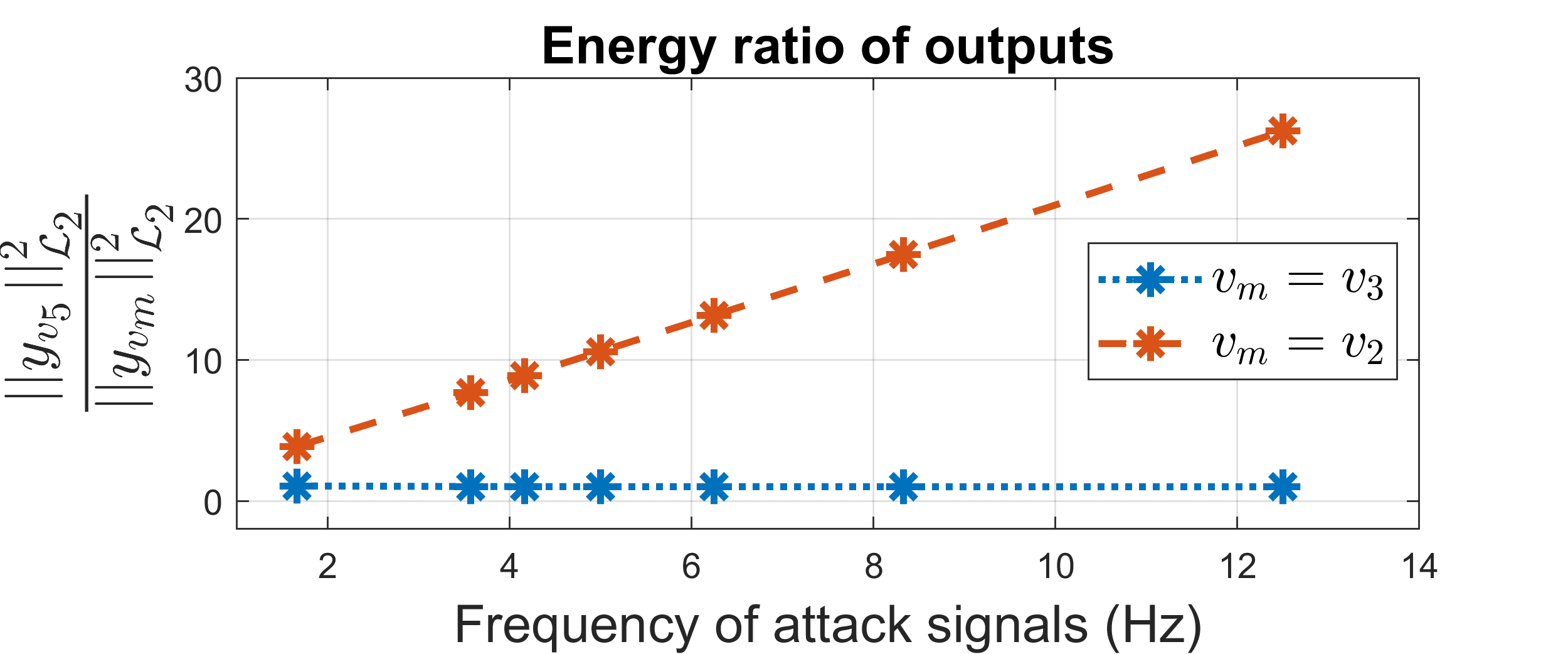}
	\caption{Energy ratio of the outputs of target $v_5$ vs. feasible monitor $v_3$ and target $v_5$ vs. infeasible monitor $v_2$.}
	\label{fig:compare}
\end{figure}

To validate the obtained results, 
let us take an example of a 10-vertex networked control system depicted in Fig. \ref{fig:graph_10vertex}. 
We simply verify that no pair of attack and monitor vertices exhibits finite unstable zeros.
Suppose that $v_5$ is  the target vertex. 
There are two feasible monitor vertices $v_3$ and $v_6$, which satisfy the algebraic condition \eqref{opt2:cond}.
Indeed, we simulate two scenarios, in which
the detector monitors the outputs of the vertices $v_3$ (feasible) and $v_2$ (infeasible).
Meanwhile,
the adversary selects the vertex $v_4$ to conduct malicious attack signals at frequency $1.67$ Hz depicted in Fig. \ref{fig:attack}. 
The outputs of the monitor vertices $v_2$, $v_3$ and the target vertex $v_5$ are shown in Fig. \ref{fig:output}.
Fig. \ref{fig:output} shows that the output of the feasible monitor vertex $v_3$ (red dash-dotted line) approximately tracks the figure for the target vertex $v_5$ (blue dashed line).
The energy produced by the output of the vertices $v_5$ and $v_3$ witnesses no noticeable difference, namely around $1.04$.
By contrast, the output energy of the infeasible monitor vertex $v_2$ (yellow dotted line) is only $0.27$, almost four times as low as the output energy of the target vertex $v_5$. 
More specifically, the output energies of those vertices over time horizon are illustrated in Fig. \ref{fig:output_energy}.
Next, we will investigate how the ratio of the output energy of the above vertices progresses when increasing the frequency of the attack signals (see Fig. \ref{fig:compare}).
%%
%The ratio of the output energy of the aforementioned vertices is given in Fig. \ref{fig:compare}.
%%
As seen  in Fig. \ref{fig:compare}, the gap between the two lines dramatically increases following the rise of the attack signal frequency.
While the blue dotted-line ($v_m=v_3$) almost remains unchanged at $1.0$, the red dashed-line ($v_m = v_2$) significantly becomes unbounded as the attack signal frequency increases.
This implies that with massively high frequencies of the attack signals, the adversary is capable of manipulating the adversarial effect on the output of the target vertex at wish  while remaining undetected at the infeasible monitor vertex $v_2$.
\begin{table*}[!t]
\caption{Game payoff \eqref{opt11} w.r.t. target vertex $v_5$ for the detector and the adversary corresponding to their chosen pair of monitor and attack vertices. \label{tab:gamma2}} 
	\centering
		%\scalebox{1}{
			\begin{tabular}{|c|c|c|c|c|c|c|c|c|c|c|}
				\hline
				\diagbox{$v_a$}{$v_m$} & $v_1$ & $v_2$ & $v_3$ & $v_4$ & $v_6$ & $v_7$ & $v_8$ & $v_9$ & $v_{10}$  & $p^\star(v_a)$ ($\%$)
				\\ \hline
				$v_1$ & 1 & 1.0062 & 1.2405 & $\infty$  & 1.4384 & $\infty$ & $\infty$ & 1.2417 & 1.0074  & 0
				\\ \hline
				$v_2$ & 1.2124 & 1 & 1.7737 & $\infty$  & 1.4669 & $\infty$ & 1.1984 & $\infty$ & $\infty$  & $\approx$ 100
				\\ \hline
				$v_3$ & 1.0565 & 1.2329 & 1 & 1.0043  & 1.1905 & 1.008 & $\infty$ & 1.2369 & 1.2681  &  $\approx$ 0
				\\ \hline
				$v_4$ & $\infty$ & $\infty$ & 1.1742 & 1  & 1.4407 & 1 & $\infty$ & 1.0126 & $\infty$ & 0 
				\\ \hline
				$v_6$ & 1.1886 & 1.0029 & 1.1729 & 1.2122 & 1 & 1.0045 & 1.212 & $\infty$ & 1.0038  & 0
				\\ \hline
				$v_7$ & $\infty$ & 2.2853 & 1 & 1  & 1 & 1 & $\infty$ & $\infty$ & 2.405  & 0
				\\ \hline
				$v_8$ & 1 & 1 & 1 & 1.2928 &  1 & $\infty$ & 1 & 1 & 1  & 0
				\\ \hline
				$v_9$ & 1 & 1 & 1 & 1 &  1 & $\infty$ & 1 & 1 & 1  & 0
				\\ \hline
				$v_{10}$ & 1 & 1 & 1 & $\infty$ &  1 & 2.3027 & 1 & 1 & 1 & 0 
				\\ \hline
				$p^\star(v_m)$ ($\%$) & 0 & $\approx$ 0 & 0 & 0 & $\approx$ 100  & 0 & 0 & 0 & 0 &
				\\ \hline
			\end{tabular}
		%}
	%\end{center}
\end{table*}

Next, the above results will be verified once again by computing the game payoff $J_5(v_a,v_m)$ with pairs of attack and monitor vertices $(v_a,v_m \neq v_5 \in \Vc)$ (see Tab. \ref{tab:gamma2}).
Looking at the third column ($v_m = v_3$) and the fifth column ($v_m = v_6$) of Tab. \ref{tab:gamma2}, no cell gives infinite value.
On the other hand, the other columns show at least one infinite game payoff.
This assessment once again confirms that  $v_m = v_3$ and $v_m = v_6$ are the feasible monitor vertices solving Problem ~\ref{prob:fea}.
%, which satisfy Def. \ref{def:set_monitor}.
%%
\\
By observing the game payoffs with the target vertex $v_5$ in Tab. \ref{tab:gamma2}, there is no pure Nash equilibrium.
However, this game always admits a mixed-strategy Nash equilibrium \citep{zhu2015game}.
Next, we investigate the mixed-strategy Nash equilibrium for this example.
Let us denote $p(v_a)$ and $p(v_m)$ as the probabilities of attack $v_a$ and monitor vertices $v_m$, respectively.
$P(v_a) = \left[p(v_{a=1}),\ldots,p(v_{a=10})\right]^\top$,
$P(v_m) = \left[p(v_{m=1}),\ldots,p(v_{m=10})\right]^\top$.
The expected game payoff w.r.t. the target vertex $v_5$ for attack vertex $v_a$ and monitor vertex $v_m$ is given by
\begin{align}
	Q_5(v_a,v_m) = P(v_a)^\top J_5 P(v_m), 
\end{align}
%%
%where $\bar{J}_5(v_a,v_m)$ is a $10\times10$-game matrix computed in Tab. \ref{tab:gamma2}.
where $J_5 = \big[J_5(v_i,v_j)_{ij}\big]$ is a $9\times9$-game matrix computed in Tab. \ref{tab:gamma2}.
There exits a saddle point $(v_a^\star,v_m^\star)$ satisfies
\begin{align*}
	-\infty < Q_5(v_a,v_m^\star) \leq Q_5(v_a^\star,v_m^\star) &\leq Q_5(v_a^\star,v_m) < \infty,
	\label{exp_gamepatoff}
	\\ 
	& \forall v_a,v_m \neq v_\tau \in \Vc.
	\non
\end{align*}
The saddle point $(v_a^\star,v_m^\star)$ in the condition above indicates that a deviation of selecting $v_a(v_m)$ does not increase(decrease) the optimal expected game payoff $Q_5(v_a^\star,v_m^\star)$. 
%For the optimal attack vertex $v_a^\star$
%%
% \begin{align}
% 	\underset{p(v_a)}{\max}  
% 	~~
% 	\underset{v_m \neq v_\tau \in \Vc}{\min} & ~~
% 	\sum_{v_a \neq v_\tau \in \Vc} J_5(v_a,v_m) ~ p(v_a)
% 	\label{minmax_va}
% 	\\
% 	\text{s.t.} &
% 	\sum_{v_a \neq v_\tau \in \Vc} p(v_a) = 1,
% 	\non \\
% 	& p(v_a) \geq 0, ~~ \forall v_a \neq v_\tau \in \Vc.
% 	\non
% \end{align}
% %%
% By using a slack variable $z$ \citep[Chapter 4]{boyd2004convex}, \eqref{minmax_va} is equivalent to
% \begin{align}
% 	\underset{p(v_a)}{\max} &~~
% 	z 
% 	\label{minmax_va1}
% 	\\
% 	\text{s.t.} & ~~
% 	{\bf 1}z - \sum_{v_a \neq v_\tau \in \Vc} J_5 ~  P(v_a) \leq 0,
% 	\non \\
% 	&\sum_{v_a \neq v_\tau \in \Vc} p(v_a) = 1,
% 	\non \\
% 	& p(v_a) \geq 0, ~~ \forall v_a \neq v_\tau \in \Vc.
% 	\non
% \end{align}
%%
% \TN{
% From Eqs. 37-38, do we need omit this part?
% }
% The same procedure in \eqref{minmax_va} and \eqref{minmax_va1}, which are used to find $P^\star(v_a)$, can also be applied to find $P^\star(v_m)$.
From the numerical results in Tab. \ref{tab:gamma2}, while the probability of selecting $v_m = 6$ is approximately 100\%, the figures for $v_a = 2$ is approximately 100\%.
The optimal probabilities $P^\star(v_a)$ and $P^\star(v_m)$ give us the optimal expected game payoff $Q_5(v_a^\star,v_m^\star) = 1.4669$.
%%
%%
%\newpage

\section{Conclusion}
\label{sec:concl}
%\tcr{[shorten conclusion to save space?]}
In this paper, we investigated a continuous-time networked control system in the presence of a cyber-attack conducted by an adversary. 
%%
%The purpose of the adversary is to manipulate the output of a protected target vertex by directly mounting the stealthy data injection attack on another vertex.
%%
%Meanwhile, 
An optimal sensor placement problem was raised such that
a detector places a sensor at a vertex to monitor such a cyber-attack. 
We invoked a single-adversary-single-detector zero-sum game to describe the optimal sensor placement problem.
%We invoked the game-theoretic approach to describe a zero-sum game between two players, the detector and the adversary.
%%
This game was then formulated by employing a min-max optimization problem.
In order to guarantee the feasibility of the min-max optimization problem,
this paper presented a necessary and sufficient condition and an  algebraic sufficient condition to find feasible monitor vertices.
By placing a sensor at one of the feasible monitor vertices, the detector possibly monitors the cyber-attack.
Further, the mixed-strategy Nash equilibrium of the zero-sum game was also analyzed to determine the optimal sensor placement.
In future works, by inheriting the concept of an untouchable target vertex in this study, our game will be expanded to consider multiple adversaries and multiple detectors.
%\begin{ack}
%Place acknowledgments here.
%\end{ack}

%\bibliographystyle{ieeetr} 
\bibliography{mybibfile}             % bib file to produce the bibliography

\begin{thebibliography}{13}
\providecommand{\natexlab}[1]{#1}
\providecommand{\url}[1]{\texttt{#1}}
\providecommand{\urlprefix}{URL }
\expandafter\ifx\csname urlstyle\endcsname\relax
  \providecommand{\doi}[1]{doi:\discretionary{}{}{}#1}\else
  \providecommand{\doi}{doi:\discretionary{}{}{}\begingroup
  \urlstyle{rm}\Url}\fi

\bibitem[{Briegel et~al.(2011)Briegel, Zelazo, B{\"u}rger, and
  Allg{\"o}wer}]{briegel2011zeros}
Briegel, B., Zelazo, D., B{\"u}rger, M., and Allg{\"o}wer, F. (2011).
\newblock On the zeros of consensus networks.
\newblock In \emph{2011 50th IEEE Conference on Decision and Control and
  European Control Conference}, 1890--1895. IEEE.

\bibitem[{Falliere et~al.(2011)Falliere, Murchu, and Chien}]{falliere2011w32}
Falliere, N., Murchu, L.O., and Chien, E. (2011).
\newblock W32. stuxnet dossier.
\newblock \emph{White paper, Symantec Corp., Security Response}, 5(6), 29.

\bibitem[{Franklin et~al.(2002)Franklin, Powell, Emami-Naeini, and
  Powell}]{franklin2002feedback}
Franklin, G.F., Powell, J.D., Emami-Naeini, A., and Powell, J.D. (2002).
\newblock \emph{Feedback control of dynamic systems}, volume~4.
\newblock Prentice hall Upper Saddle River, NJ.

\bibitem[{Gupta et~al.(2016)Gupta, Langbort, and
  Ba{\c{s}}ar}]{gupta2016dynamic}
Gupta, A., Langbort, C., and Ba{\c{s}}ar, T. (2016).
\newblock Dynamic games with asymmetric information and resource constrained
  players with applications to security of cyberphysical systems.
\newblock \emph{IEEE Transactions on Control of Network Systems}, 4(1), 71--81.

\bibitem[{Khalil(2002)}]{khalil2002nonlinear}
Khalil, H.K. (2002).
\newblock Nonlinear systems third edition.
\newblock \emph{Patience Hall}, 115.

\bibitem[{Miao et~al.(2018)Miao, Zhu, Pajic, and Pappas}]{miao2018hybrid}
Miao, F., Zhu, Q., Pajic, M., and Pappas, G.J. (2018).
\newblock A hybrid stochastic game for secure control of cyber-physical
  systems.
\newblock \emph{Automatica}, 93, 55--63.

\bibitem[{Morris and Rebarber(2010)}]{morris2010invariant}
Morris, K. and Rebarber, R. (2010).
\newblock Invariant zeros of siso infinite-dimensional systems.
\newblock \emph{International journal of control}, 83(12), 2573--2579.

\bibitem[{Pirani et~al.(2021)Pirani, Nekouei, Sandberg, and
  Johansson}]{pirani2021game}
Pirani, M., Nekouei, E., Sandberg, H., and Johansson, K.H. (2021).
\newblock A game-theoretic framework for the security-aware sensor placement
  problem in networked control systems.
\newblock \emph{IEEE Transactions on Automatic Control}.

\bibitem[{Teixeira et~al.(2015{\natexlab{a}})Teixeira, Sandberg, and
  Johansson}]{teixeira2015strategic}
Teixeira, A., Sandberg, H., and Johansson, K.H. (2015{\natexlab{a}}).
\newblock Strategic stealthy attacks: the output-to-output $\ell_2$-gain.
\newblock In \emph{2015 54th IEEE Conference on Decision and Control (CDC)},
  2582--2587. IEEE.

\bibitem[{Teixeira et~al.(2015{\natexlab{b}})Teixeira, Shames, Sandberg, and
  Johansson}]{teixeira2015secure}
Teixeira, A., Shames, I., Sandberg, H., and Johansson, K.H.
  (2015{\natexlab{b}}).
\newblock A secure control framework for resource-limited adversaries.
\newblock \emph{Automatica}, 51, 135--148.

\bibitem[{Teixeira(2021)}]{teixeira2021security}
Teixeira, A.M. (2021).
\newblock Security metrics for control systems.
\newblock In \emph{Safety, Security and Privacy for Cyber-Physical Systems},
  99--121. Springer.

\bibitem[{Torres and Roy(2015)}]{torres2015graph}
Torres, J.A. and Roy, S. (2015).
\newblock Graph-theoretic analysis of network input--output processes: Zero
  structure and its implications on remote feedback control.
\newblock \emph{Automatica}, 61, 73--79.

\bibitem[{Zhu and Basar(2015)}]{zhu2015game}
Zhu, Q. and Basar, T. (2015).
\newblock Game-theoretic methods for robustness, security, and resilience of
  cyberphysical control systems: games-in-games principle for optimal
  cross-layer resilient control systems.
\newblock \emph{IEEE Control Systems Magazine}, 35(1), 46--65.

\end{thebibliography}

\end{document}